\documentclass[english,11pt]{amsart}
\usepackage[english]{babel}
\usepackage{amsmath,amssymb,enumerate,amsthm}
\usepackage[latin1]{inputenc}
\usepackage[T1]{fontenc}
\usepackage{cite}
\usepackage{latexsym}
\usepackage{setspace}

\newtheorem{lemma}{Lemma}[section]
\newtheorem{teorema}[lemma]{Theorem}
\newtheorem{prop}[lemma]{Proposition}
\newtheorem{coro}[lemma]{Corollary}

\newtheorem{lema}[lemma]{Lemma}
\newtheorem{rk}[lemma]{Remark}

\newcommand{\SL}{\mathrm{SL}(2,\mathbb{Z})}

\def\={\;=\;}
\def\.={\;\dot{=}\;}
\newcommand{\Gal}{\mathrm{Gal}}
\def\O{\mathcal{O}}

\begin{document}

\title{Galois action on special theta values}
\author{Paloma Bengoechea}
\address{Department of Mathematics, University of York, York, YO10 5DD, United Kingdom}
\email{paloma.bengoechea@york.ac.uk}

\maketitle


\begin{abstract}
For a primitive Dirichlet character $\chi$ of conductor $N$ set $\theta_\chi(\tau)\=\sum_{n\in\mathbb{Z}}n^\epsilon\, \chi(n)\, e^{\pi i n^2\tau/N}$ (where $\epsilon = 0$ for even $\chi$, $\epsilon = 1$ for odd $\chi$) the associated theta series. Its value at its point of symmetry
under the modular transformation $\tau\mapsto -1/\tau$ is related by $\theta_\chi(i)=W(\chi)\theta_{\bar{\chi}}(i)$ to the
root number of the L-series of $\chi$ and hence can be used to calculate the latter quickly if it does not
vanish. 
Using Shimura's reciprocity law, we calculate the Galois action on these special values of \textit{theta} functions with odd $N$ normalised by the Dedekind eta function. As a consequence, we prove some experimental results of Cohen and Zagier and we deduce a partial result on the non-vanishing of these special theta values with prime $N$.
\end{abstract}

\section{Introduction}

\let\thefootnote\relax\footnote{2010 Mathematics Subject Classification: 11G15, 14K22, 14K25, 11M06.\\Key words: theta functions, complex multiplication, L-series, Shimura's reciprocity law.}

Let $\chi$ be a primitive Dirichlet character with conductor $N$ and order $m$. The $theta$ series associated to the character $\chi$ is defined on the upper half-plane $\mathcal{H}$ by
\begin{equation}
\theta_\chi(\tau)\=\sum_{n\in\mathbb{Z}} n^\epsilon\, \chi(n)\, q^{n^2/2N}\qquad(q=\underline{e}(\tau), \, \tau\in\mathcal{H}),
\end{equation}
where $\epsilon$ equals 0 if $\chi$ is even, or 1 if $\chi$ is odd. (Here and in the sequel we use the notation $\underline{e}(x)=e^{2i\pi x}$ for $x\in\mathbb{C}$.)
The theta series satisfies the functional equation
\begin{equation}\label{ec funcional}
\theta_{\bar\chi}(-1/\tau)\=W(\chi)\;(\tau/i)^{1/2+\epsilon}\;\theta_{\chi}(\tau),
\end{equation}
where $W(\chi)$ is the algebraic number of module 1 called \textsl{root number}, defined in an explicit way as the normalized Gauss sum associated to  $\chi$
\begin{equation}\label{def W}
W(\chi)\=G(\chi)/\sqrt{N},\qquad G(\chi)\=\sum_{n\, \mathrm{mod}\, N}\chi(n)\;\underline{e}(n/N).
\end{equation}
From equation \eqref{ec funcional} one deduces the analytic continuation and the functional equation of the $L$-series $L(s,\chi)=\sum_{n\in\mathbb{Z}}\chi(n)n^{-s}$. For a given value $s$, we cannot compute $L(s,\chi)$ with a big precision directly from its definition because it is very slowly convergent or even not convergent at all. However, we can use its approximative functional equation which arises from truncating the series; in this case, we can compute $L(s,\chi)$ in $O(\sqrt{N})$ time if the value of $W(\chi)$ is known. 
From definition \eqref{def W}, we compute $W(\chi)$ in $O(N)$ steps. In fact we can do better: considering equation \eqref{ec funcional} with the value $\tau\=i$, we deduce the identity
$$ 
W(\chi)\=\dfrac{\theta_{\bar\chi}(i)}{\theta_\chi(i)},
$$
from which we compute $W(\chi)$ in $O(\sqrt{N})$ steps when $\theta_\chi(i)\neq 0$. A natural question arises then: does $\theta_\chi(i)$ vanish for any $\chi$?

Louboutin proved in \cite{Lo} that there exists a constant $c>0$ such that, for every prime $p$, at least $cp/\log p$ of the $(p-1)/2$ values $\theta_\chi(i)$, where $\chi$ is odd with prime conductor $p$, do not vanish. Cohen and Zagier described explicit computational results in \cite{C-Z} showing that $\theta_\chi(i)\neq 0$ for the first 500 millions of characters $\chi$ with $N\leq 52100$, except for exactly (up to complex conjugation) two even characters with respective conductor 300 and 600. 

Moreover, they defined the functions
$$
A_\chi(\tau)=\dfrac{\theta_\chi(\tau/N)}{\eta(\tau/N)^{1+2\epsilon}},\qquad B_\chi(\tau)=|A_\chi(\tau)|^2=A_\chi(\tau)\, A_{\bar{\chi}}(\tau),
$$
where $\eta$ is the Dedekind's eta function defined by
$$
\eta(\tau)\=q^{1/24}\prod_{n\=1}^\infty(1-q^n),
$$
and studied the algebraic numbers $A_\chi(ip)$ and $B_\chi(ip)$ when $p$ is prime and it is the conductor of $\chi$. Indeed, since the functions above are modular functions, they are algebraic on the points of complex multiplication. Because of the algebraicity, the numbers $B_\chi(iN)$ are much easier to study than the values $\theta_\chi(i)$. Also is the product of the numbers $B_\chi(iN)$ for all characters $\chi$ with fixed conductor $N$ and fixed order $m$ (up to complex conjugation). We denote these products by
$$
\mathcal{N}(N,m)=\prod_{\substack{\mathrm{order}(\chi)= m\\ \chi\approx\bar\chi}} B_\chi(iN).
$$ 
Cohen and Zagier speculated that the values $\mathcal{N}(p,m)^2$ always belong to $\mathbb{Q}(i,j(ip))$. Moreover, if we denote by $\mathcal{N}(p,m)^d$ the smallest power of $\mathcal{N}(p,m)$ belonging to $\mathbb{Q}(i,j(ip))$, then the experimental results led Cohen and Zagier to conjecture that, for the special case of the trivial character, $d=1$ if $p\equiv 1\pmod{4}$ and $d=2$ if $p\equiv 3\pmod{4}$; for Legendre's character ($m=2$), it seems that $d=1$.

Concerning the numbers $A_\chi(ip)$, Cohen and Zagier observed that the degree drastically decreases for some powers. If we denote by $\zeta_m$ the $m$-th root of unity $e^{2\pi im}$ and $\sigma_s$ the element of $\mathrm{Gal}(\mathbb{Q}(i,j(ip),\zeta_m)/\mathbb{Q}(i,j(ip))$ sending $\zeta_m$ to $\zeta^s_m$, they speculated  $A_\chi(ip)^k\in\mathbb{Q}(j(ip),\zeta_m)$ for some $k\in\mathbb{N}$ and $A_{\chi^s}(ip)^k=\sigma_s(A_\chi(ip)^k)$ for all $s\in(\mathbb{Z}/m\mathbb{Z})^\ast$.

We are able to calculate the Galois action on these algebraic numbers and, using class field theory and Shimura's reciprocity law, we prove Cohen and Zagier's experimental results mentioned above and the generalizations to odd conductors. Concerning the non-vanishing of the special theta values, we prove $\theta_\chi(i)\neq 0$ 
for all non-quadratic $\chi$ with prime and ``big'' conductor $p=2l+1$, where $l$ is also prime (so $l$ is a Sophie Germain prime).

\section{Modularity}

Throughout the paper we denote by $\chi$ a primitive character with odd conductor $N$ and order $m$.

In this section we explicit the action of the group $\Gamma_\theta\cap\Gamma_0(N)$ on $\theta_\chi(\tau)$.
We can decompose the theta series in the following way:
\begin{equation}\label{decomposicion}
\theta_\chi(\tau)\=\sum_{h\, \mathrm{mod}\, N}\chi(h)\, \theta_{N,h}^{(\epsilon)}(\tau),
\end{equation}
where the coefficients $\chi(h)$ are $m$-th roots of unity and
\begin{equation}\label{theta parcial}
\theta_{N,h}^{(\epsilon)}(\tau)=\sum_{\substack{n\in\mathbb{Z} \\n\equiv h\!\pmod{N}}} n^\epsilon\, q^{n^2/2N}.
\end{equation} 

We define the group
$$
\Gamma_\theta\=\left\{\begin{pmatrix}a &b\\c &d\end{pmatrix}\in\SL\, :\, \begin{pmatrix}a &b\\c &d\end{pmatrix}\equiv\begin{pmatrix}1 &0\\0 &1\end{pmatrix}\mbox{ or }\begin{pmatrix}0 &1\\1 &0\end{pmatrix}\pmod{2}\right\}.
$$

In order to compute the action of $\Gamma_\theta\cap\Gamma_0(N)$ on the functions $\theta^{(\epsilon)}_{N,h}(\tau)$, we use Proposition 10.4 in \cite{I}, namely
\begin{prop}[Iwaniec]\label{inversion theta} We have
\begin{equation}\label{inversion}
\theta_{N,h}^{(\epsilon)}(-1/\tau)\=(i/N)^{1/2}\, (-\tau)^{1/2+\epsilon}\, \sum_{l\, \mathrm{mod}\, N}\underline{e}(hl/N)\, \theta_{N,l}^{(\epsilon)}(\tau).
\end{equation}
\end{prop}

\begin{prop}\label{Iwanieck} For $\gamma=\begin{pmatrix}a &b\\c &d\end{pmatrix}\in\Gamma_\theta\cap\Gamma_0(N)$, we have 
\begin{equation}\label{ec theta}
\theta^{(\epsilon)}_{N,h}(\gamma(\tau))\=\underline{e}\Big(\dfrac{a^2bdh^2}{2N}\Big)\, \upsilon(\gamma,N)\, (c\tau+d)^{1/2+\epsilon}\,  \theta^{(\epsilon)}_{N,ah}(\tau),
\end{equation}
with
\begin{equation}\label{zeta}
\upsilon(\gamma,N)\=\left\{\begin{array}{ll}\zeta_8^{bN}\left(\dfrac{d}{|bN|}\right)&\mbox{if $d$ is even},\\
\zeta_8^{d-1}\left(\dfrac{-bN}{d}\right)&\mbox{if $d$ is odd},
\end{array}\right.
\end{equation}
where $\left(\frac{\cdot}{\cdot}\right)$ is the Kronecker symbol.
\end{prop}

\textbf{Proof.}  First we suppose $d>0.$ We write
$$
\gamma'=\begin{pmatrix}a &b\\c &d\end{pmatrix}\begin{pmatrix}0 &-1\\1 &0\end{pmatrix}=\begin{pmatrix}b &-a\\d &-c\end{pmatrix}.
$$
Since $d\gamma'(\tau)=b-\dfrac{1}{d\tau-c}$, we have
\begin{equation}\label{theta eps}
\theta^{(\epsilon)}_{N,h}(\gamma'(\tau))\=\sum_{n\equiv h\pmod{N}} n^\epsilon\, \underline{e}\Big(\dfrac{n^2}{2N}\Big(\dfrac{b}{d}-\dfrac{1}{d(d\tau-c)}\Big)\Big).
\end{equation}
But $\underline{e}(bn^2/2dN)$ only depends on $n$ (mod $dN$). Indeed, if $b\equiv 0\pmod{2}$, then this assertion is obvious. Otherwise, $d\equiv 0\pmod{2}$ and for $n=kdN+r$ with $1\leq r\leq dN$, $k\in\mathbb{Z}$, we have $n^2\equiv r^2\pmod{2dN}$.

Hence we can split the sum \eqref{theta eps} into classes modulo $dN$:
$$
\theta^{(\epsilon)}_{N,h}(\gamma'(\tau))\=\sum_{\substack{m\, \mathrm{mod}\, dN\\ m\equiv h\pmod{N}}}\underline{e}\Big(\dfrac{bm^2}{2dN}\Big)\, \sum_{n\equiv m\pmod{dN}}n^\epsilon\, \underline{e}\Big(\dfrac{n^2}{2dN}\dfrac{-1}{d\tau-c}\Big).
$$
The second sum is the theta function associated (in the sense \eqref{theta parcial}) to the conductor $dN$ and residual class $m$ (mod $dN$) evaluated on $\dfrac{-1}{d\tau-c}$. By applying Proposition \ref{inversion theta} to this sum, we obtain
$$
\Big(\dfrac{i}{dN}\Big)^{1/2}\, (c-d\tau)^{1/2+\epsilon}\, \sum_{l\, \mathrm{mod}\, dN}\underline{e}\Big(\dfrac{lm}{dN}\Big)\, \sum_{n\equiv l\pmod{dN}}n^\epsilon\, \underline{e}\Big(\dfrac{n^2}{2dN}(d\tau-c)\Big).
$$
If $d\equiv 0\pmod{2}$, then $n^2\equiv l^2\pmod{2dN}$. Otherwise, $c\equiv 0\pmod{2}$. In both situations, $cn^2\equiv cl^2\pmod{2dN}$. Thus
\begin{equation}\label{3}
\theta^{(\epsilon)}_{N,h}(\gamma'(\tau))=\Big(\dfrac{i}{dN}\Big)^{1/2} (c-d\tau)^{1/2+\epsilon} \sum_{l\, \mathrm{mod}\, dN}\varphi(h,l)\, \sum_{n\equiv l\!\!\pmod{dN}}n^\epsilon\, \underline{e}\Big(\dfrac{n^2}{2N}\tau\Big),
\end{equation}
where
$$
\varphi(h,l)=\sum_{\substack{m\, \mathrm{mod}\, dN\\ m\equiv h\pmod{N}}}\underline{e}((bm^2+2lm-cl^2)/2dN).
$$
We rewrite $\varphi(h,l)$ after changing the variable $m$ by $m+cl$:
\begin{align*}
\varphi(h,l)&\=\sum_{\substack{m\, \mathrm{mod}\, dN\\ m\equiv h-cl\pmod{N}}}\underline{e}((b(m+cl)^2+2l(m+cl)-cl^2)/2dN)\\
&\=\sum_{\substack{m\, \mathrm{mod}\, dN\\ m\equiv h-cl\pmod{N}}}\underline{e}((bm^2+2adlm+acdl^2)/2dN)
\end{align*}
since $ad-bc=1$.
In the term $2adlm$, we replace $m$ by $h-cl$ (mod $N$):
\begin{equation}\label{phi}
\varphi(h,l)\=\underline{e}(2ahl-acl^2/2N)\, \varphi(h-cl,0).
\end{equation}
This expression makes possible to replace in \eqref{3} $l$ (mod $dN$) by $l$ (mod $N$). We obtain 
\begin{equation}\label{ecuacion final}
\theta^{(\epsilon)}_{N,h}(\gamma'(\tau))\=\Big(\dfrac{i}{dN}\Big)^{1/2}\, (c-d\tau)^{1/2+\epsilon}\, \sum_{l\, \mathrm{mod}\, N}\varphi(h,l)\, \theta^{(\epsilon)}_{N,l}(\tau).
\end{equation}

Replacing $\tau$ by $-1/\tau$ and applying Proposition \ref{inversion theta} to each $\theta^{(\epsilon)}_{N,l}(-1/\tau)$, we get
$$
\theta^{(\epsilon)}_{N,h}(\gamma(\tau))\=\dfrac{(-1)^\epsilon}{d^{1/2}N}\, (c\tau+d)^{1/2+\epsilon}\, \sum_{l\, \mathrm{mod}\, N}\phi(h,l)\, \theta^{(\epsilon)}_{N,l}(\tau), 
$$
where
$$
\phi(h,l)\=\sum_{g\, \mathrm{mod}\, N}\varphi(h,g)\, \underline{e}(gl/N).
$$

Since $c\equiv 0\pmod{N}$ and $ac\equiv 0\pmod{2}$, the formula \eqref{phi} becomes
$$
\varphi(h,l)=\underline{e}(ahl/N)\, \varphi(h,0).
$$
Hence
\begin{align*}
\phi(h,l)&\=\varphi(h,0)\sum_{g\, \mathrm{mod}\, N}\underline{e}(g(ah+l)/N)\\
&\=\left\{\begin{array}{ll}\varphi(h,0)N &\mbox{if $l\equiv -ah\pmod{N}$}\\ 0 &\mbox{otherwise}.
\end{array}\right.
\end{align*}

Therefore
\begin{align*}
\theta^{(\epsilon)}_{N,h}(\gamma(\tau))&\=\dfrac{(-1)^\epsilon}{d^{1/2}}\, \varphi(h,0)\, (c\tau+d)^{1/2+\epsilon}\,  \theta^{(\epsilon)}_{N,-ah}(\tau)\\
&\\
&\=\dfrac{\varphi(h,0)}{d^{1/2}}\, (c\tau+d)^{1/2+\epsilon}\, \theta^{(\epsilon)}_{N,ah}(\tau).
\end{align*}

We still have to calculate
$$
\varphi(h,0)=\sum_{\substack{m\, \mathrm{mod}\, dN\\m\equiv h\pmod{N}}}\underline{e}\Big(\dfrac{bm^2}{2dN}\Big).
$$
Since $ad\equiv 1\pmod{N}$, we can write $m=adh+nN$ with $1\leq n\leq d$. Thus we get
$$
\varphi(h,0)=\underline{e}\Big(\dfrac{a^2bdh^2}{2N}\Big)S_{bN,d},
$$
where
$$
S_{bN,d}=\sum_{1\leq n\leq d}\underline{e}\Big(\dfrac{bNn^2}{2d}\Big)
$$
is a well known Gauss sum, calculated for example in \cite{Mu}: 
$$
S_{bN,d}=\left\{\begin{array}{ll} d^{1/2}\, \zeta_8^{bN}\left(\dfrac{d}{|bN|}\right) &\quad\mbox{if $d$ is even},\\
&\\
d^{1/2}\, \zeta_8^{d-1}\left(\dfrac{-bN}{d}\right) &\quad\mbox{if $d$ is odd}.
\end{array}\right.
$$
Finally we obtain \eqref{ec theta} for $d>0$. When $d<0$, we can change $\gamma$ by $-\gamma$ such that the left-hand term of the equality \eqref{ec theta} does not vary. It is easily shown that the right-hand term does not vary either, i.e, $\upsilon(-\gamma,N)i=\upsilon(\gamma,N)$.
\begin{flushright}
$\square$
\end{flushright}

Meyer's formula (\cite{Me}) gives, for $\gamma=\begin{pmatrix}a &b\\c &d\end{pmatrix}\in\SL$, some functions $\epsilon_1(\gamma)$ and $\epsilon_2(\gamma)$ such that
\begin{equation}\label{meyer}
\eta(\gamma(\tau))=\epsilon_1(\gamma)\, \epsilon_2(\gamma)\, (c\tau+d)^{1/2}\, \eta(\tau).
\end{equation}
We can fix $c>0$ or $c=0$ and $d>0$, changing $\gamma$ by $-\gamma$ if necessary; then $\mathrm{Im}(c\tau+d)\geq 0$ and we chose $\mathrm{Re}(c\tau+d)^\frac{1}{2}\geq 0$. If $c>0$, we write $c=2^r\cdot c_0$ with $c_0$ odd. If $c=0$, we write $c_0=r=1$. Then we have
$$
\epsilon_1(\gamma)=\left(\dfrac{a}{c_0}\right)\quad\mbox{and}\quad \epsilon_2(\gamma)=\zeta_{24}^{ab+cd(1-a^2)-ca+3c_0(a-1)+r\frac{3}{2}(a^2-1)}.
$$

\begin{prop}\label{coro inv theta}
Let $w=\dfrac{24N}{(12,N)}$. 
The functions $\dfrac{\theta^{(\epsilon)}_{N,h}(\tau)}{\eta^{1+2\epsilon}(\tau)}$ are $\Gamma(w)$-invariant and the functions $\dfrac{\theta^{(\epsilon)}_{N,h}(\tau/N)}{\eta^{1+2\epsilon}(\tau/N)}$ are $\Gamma(wN)$-invariant. 
\end{prop}

\textbf{Proof.} For $\gamma=\begin{pmatrix}a &b\\c &d\end{pmatrix}\in\Gamma(w)$, the multiplicative system $\upsilon(\gamma,N)$ in Proposition \ref{Iwanieck} becomes simpler (see \cite{I} Proposition 10.6): $\upsilon(\gamma,N)=\epsilon_1(\gamma)$. The same happens with the second Meyer's function: $\epsilon_2(\gamma)=1$. Hence the functions $\theta^{(\epsilon)}_{N,h}(\tau)/\eta^{1+2\epsilon}(\tau)$ are $\Gamma(w)$-invariant.

Let $\gamma=\begin{pmatrix}a &b\\c &d\end{pmatrix}$ be an element in $\Gamma(wN)$. We write
$$
\gamma'\={\begin{pmatrix}1 &0\\0 &N\end{pmatrix}\begin{pmatrix}a &b\\c &d\end{pmatrix}\begin{pmatrix}1 &0\\0 &N\end{pmatrix}^{-1}\=\begin{pmatrix}a &\frac{b}{N}\\cN &d\end{pmatrix}},
$$
such that $\gamma'\in\Gamma(w)$ and
$$
\dfrac{\theta^{(\epsilon)}_{N,h}(\gamma(\tau)/N)}{\eta^{1+2\epsilon}(\gamma(\tau)/N)}=\dfrac{\theta^{(\epsilon)}_{N,h}(\gamma'(\tau/N))}{\eta^{1+2\epsilon}(\gamma'(\tau/N))}=\dfrac{\theta^{(\epsilon)}_{N,h}(\tau/N)}{\eta^{1+2\epsilon}(\tau/N)}.
$$
\begin{flushright}
$\square$
\end{flushright}

\section{Shimura's reciprocity law}

In this section we follow the interpretation of Shimura's reciprocity law (see \cite{Sh}) by Gee and Stevenhagen (see \cite{GS}, \cite{G1}, \cite{St}).
Let $K$ be an imaginary quadratic field and $\O$ an order in $K$ with basis $[\alpha,1]$. The first fundamental theorem of complex multiplication states that the $j$-invariant $j(\alpha)$ is an algebraic integer and $K(j(\alpha))$ is the ring class field $H_\mathcal{O}$ of $\mathcal{O}$ (see, for example, \cite{C}).
For $M\geq 1$, the field $F_M$ of modular functions with level $M$ is defined as the field of meromorphic functions on $\mathcal{H}\cup\left\{\infty\right\}$, invariant by $\Gamma(M)$ and whose coefficients in the Fourier expansion in the variable $q^{1/M}$ belong to the field $\mathbb{Q}(\zeta_M)$. It follows from the second fundamental theorem of complex multiplication, stated for example in \cite{C} and proved in \cite{La2} and \cite{Fra}, that for a function $f$ belonging to the field $F_M$, the value $f(\alpha)$ is an element of the ray class field $H_{M,\mathcal{O}}$ with conductor $M$ over the ring class field $H_\mathcal{O}$. 

Shimura's reciprocity law gives the action of the group $\Gal(H_{M,\O}/H_\O)$ on $f(\alpha)$ combining Artin's reciprocity law arisen from class field theory, and Galois theory on $F_M$. Artin's reciprocity law gives the exact sequence
$$
 \O^\ast\longrightarrow(\O/M\O)^\ast\stackrel{A}{\longrightarrow}\Gal(H_{M,\O}/H_\O)\longrightarrow 1,
$$
where $A$ is the Artin map.
The map
$$
\begin{array}{rll}
\mathrm{GL}_2(\mathbb{Z}/M\mathbb{Z})&\longrightarrow&\Gal(F_M/F_1)\\
\mu=\tiny{\begin{pmatrix}1 &0\\0 &\det(\mu)\end{pmatrix}}\gamma&\mapsto &(\sum c_k q^{k/M}\mapsto (\sum \sigma_{\det(\mu)}(c_k)q^{k/M})|_0\gamma),
\end{array}
$$
where $\gamma\in\mathrm{SL}_2(\mathbb{Z}/M\mathbb{Z})$ and $\sigma_{\det(\mu)}\in\mathrm{Aut}(\mathbb{Q}(\zeta_M))$ sends $\zeta_M$ to $\zeta_M^{\det(\mu)}$, is surjective. When $D<-4$, its kernel is $\left\{\pm 1\right\}$.

Then we have the following diagram, where all the sequences are exact:
$$
\begin{array}{ccccccc}
 \O^\ast&\longrightarrow&(\O/M\O)^\ast&\stackrel{A}{\longrightarrow}&\Gal(H_{M,\O}/H_\O)&\longrightarrow &1\\
 &&\downarrow g_\alpha\\
 \left\{\pm 1\right\}&\longrightarrow&\mathrm{GL}_2(\mathbb{Z}/M\mathbb{Z})&\longrightarrow&\Gal(F_M/F_1)&\longrightarrow&1.
\end{array} 
$$
The connection map $g_\alpha$ sends $x\in(\O/M\O)^\ast$ to the matrix corresponding to the multiplication by $x$ with respect to the basis $[\alpha,1]$ ($g_\alpha(x)\tiny{\begin{pmatrix}\alpha\\1\end{pmatrix}}=\tiny{\begin{pmatrix}x\alpha\\x\end{pmatrix}}$). If $X^2+Bx+C$ is the irreducible polynomial of $\alpha$ over $\mathbb{Q}$, we can explicitely describe $g_\alpha$ by
$$
\begin{array}{rrll}
g_\alpha:&(\O/M\O)^\ast&\longrightarrow &GL_2(\mathbb{Z}/M\mathbb{Z})\\
&x=s\alpha+t &\mapsto &\begin{pmatrix}t-Bs &-Cs\\s &t\end{pmatrix}.
\end{array}
$$
The map $g_\alpha$ gives an action of $(\O/M\O)^\ast$ on $F_M$ and the reciprocity relation: for $x\in(\mathcal{O}/M\mathcal{O})^\ast$,
$$
(f(\alpha))^x=(f^{g_\alpha(x^{-1})})(\alpha).
$$
Moreover, denoting by $F=\bigcup_{M\geq 1} F_M$ the modular field, if the extension $F/\mathbb{Q}(f)$ is Galois, then
we have the fundamental equivalence:
$$
(f(\alpha))^x=f(\alpha)\quad \Leftrightarrow\quad f^{g_\alpha(x)}=f.
$$
We denote by
$$
W_{M,\alpha}\=\left\{\begin{pmatrix}t-Bs &-Cs\\s &t\end{pmatrix}\in\mathrm{GL}_2(\mathbb{Z}/M\mathbb{Z})\mid t,s\in\mathbb{Z}/M\mathbb{Z}\right\}
$$
the image of $(\O/M\O)^\ast$ by $g_\alpha$ when $D<-4$. The algebraic number $f(\alpha)$ belongs to $H_\O$ if $f$ is invariant by the action of $W_{M,\alpha}/\left\{\pm 1\right\}$.

\section{Galois action, proofs of the experimental results}

Let $\chi$ be a primitive character with odd conductor $N$ and order $m$. By Proposition \ref{coro inv theta}, the functions $\frac{\theta_{N,h}^{(\epsilon)}(\tau)}{{\eta^{1+2\epsilon}}(\tau)}$ belong to the field $F_w$, where $w=\frac{24}{(12,N)}$. Hence we deduce (see decomposition \eqref{decomposicion}) that the numbers $A_{\chi}(iN)$ belong to the field $H_{w,\O_K}(\zeta_m)$, where $\O_K$ is the ring of integers of the field $K=\mathbb{Q}(i)$. In this section we use Shimura's reciprocity law to obtain more accurated statements about the algebraicity of the numbers $A_\chi(iN)$ and $B_\chi(iN)$. 

Let
$$
v=\chi(-1),\qquad M=24mN^2,
$$
and $n=m$ if $m$ is even and $2m$ otherwise. We consider the order $\O=\mathbb{Z}[iN]$ in $K=\mathbb{Q}(i)$, and its ring class field $H_\O=K(j(iN))$.

By Proposition \ref{coro inv theta}, we know that the functions $A_\chi(\tau)$ and $B_{\chi}(\tau)$ belong to the field $F_{M}$. Following the notations of section 3, $$
W_{M,iN}=\left\{\begin{pmatrix} t &-N^2s\\s &t\end{pmatrix}\in\mathrm{GL}_2(\mathbb{Z}/M\mathbb{Z})\, \mid\, t,s\in\mathbb{Z}/M\mathbb{Z}\right\}.
$$ 

\begin{prop}\label{prop 2} For $\mu=\begin{pmatrix}t &-N^2s\\s &t\end{pmatrix}\in W_{M,iN}$, 
we have
$$
(B_{\chi}|\mu)(\tau)\=(-1)^{\frac{N-v}{2}(t-1)}\, B_{\chi^{\mathrm{det}(\mu)}}(\tau)
$$
and
$$
(A_{\chi}|\mu)(\tau)^n\=(-1)^{\frac{(N-v)n}{2}(t-1)}\, A_{\chi^{\mathrm{det}(\mu)}}(\tau)^n.
$$
\end{prop}

\textbf{Proof.} Let $\mu=\begin{pmatrix}t &-N^2s\\s &t\end{pmatrix}$ be an element in $W_{M,iN}$.
We write
$$
\mu\={\begin{pmatrix}1 &0\\0 &\mathrm{det}(\mu)\end{pmatrix}\begin{pmatrix}t &-N^2s\\s(\mathrm{det}(\mu))^{-1} &t(\mathrm{det}(\mu))^{-1}\end{pmatrix}}.
$$
The first matrix transforms $B_\chi(\tau)$ into $B_{\chi^{\mathrm{det}(\mu)}}(\tau)$. To explicit the action of the second matrix we chose  
 $\gamma\=\begin{pmatrix}a &b\\c &d\end{pmatrix}\in \SL$ a representant of $\begin{pmatrix}t &-N^2s\\s(\mathrm{det}(\mu))^{-1} &t(\mathrm{det}(\mu))^{-1}\end{pmatrix}\in\mathrm{SL}_2(\mathbb{Z}/M\mathbb{Z})$ with $c>0$, or $c=0$ and $d>0$.
Since
\begin{equation}\label{congruencias 1}
a\equiv d\, \mathrm{det}(\mu)\pmod{M},\qquad b\equiv -cN^2\, \mathrm{det}(\mu)\pmod{M}
\end{equation}
and $N$ is odd, we have $\gamma\in\Gamma_\theta\cap\Gamma^0(N^2)$. 

We write
\begin{equation}\label{a}
\gamma'\={\begin{pmatrix}1 &0\\0 &N\end{pmatrix}\begin{pmatrix}a &b\\c &d\end{pmatrix}\begin{pmatrix}1 &0\\0 &N\end{pmatrix}^{-1}\=\begin{pmatrix}a &\frac{b}{N}\\cN &d\end{pmatrix}}
\end{equation}
such that $\gamma'$ satisfies the conditions of Proposition \ref{Iwanieck} and 
\begin{equation}\label{B}
B_{\chi}(\gamma(\tau))=\dfrac{\theta_{\chi}\Big(\gamma'\Big(\dfrac{\tau}{N}\Big)\Big)\theta_{\bar{\chi}}\Big(\gamma'\Big(\dfrac{\tau}{N}\Big)\Big)}{\eta\Big(\gamma'\Big(\dfrac{\tau}{N}\Big)\Big)^{2(1+2\epsilon)}}.
\end{equation}

Meyer's formula \eqref{meyer} gives
$$
\eta\Big(\gamma'\Big(\frac{\tau}{N}\Big)\Big)^2\=\epsilon_1(\gamma')^2\, \epsilon_2(\gamma')^2\, (c\tau+d)\, \eta\Big(\frac{\tau}{N}\Big)^2
$$
with $\epsilon_1(\gamma')^2=1$ and $\epsilon_2(\gamma')^2=\zeta_{12}^{a\frac{b}{N}+cdN(1-a^2)-acN+3c_0N(a-1)}$, where $c=2^r c_0$ with $c_0$ odd if $c>0$, and $c_0=1$ if $c=0$.

On the other hand, Proposition \ref{Iwanieck} gives the expression for the numerator of \eqref{B}:
\begin{equation}\label{expresion}
\upsilon(\gamma',N)^2\, (c\tau+d)^{2+4\epsilon}\sum_{h_1, h_2\, \mathrm{mod}\, N} \chi(a^{-1}h_1)\, \bar\chi(a^{-1}h_2)\, \theta_{N,h_1}^{(\epsilon)}\Big(\dfrac{\tau}{N}\Big)\, \theta_{N,h_2}^{(\epsilon)}\Big(\dfrac{\tau}{N}\Big),
\end{equation}
where
$$
\upsilon(\gamma',N)=\upsilon(\gamma,1)=\left\{\begin{array}{ll}
\zeta^b_8\, \left(\dfrac{d}{|b|}\right) &\quad\mbox{if $d$ is even}\\
\zeta^{d-1}_8\, \left(\dfrac{-b}{d}\right) &\quad\mbox{if $d$ is odd}.
\end{array}\right.
$$
Since $\chi(t^{-1}h_1)\bar\chi(t^{-1}h_2)=\chi(h_1)\bar\chi(h_2)$, the numerator of \eqref{B} becomes
$$
\theta_\chi\Big(\gamma'\Big(\dfrac{\tau}{N}\Big)\Big)\, \theta_{\bar\chi}\Big(\gamma'\Big(\dfrac{\tau}{N}\Big)\Big)\=\upsilon(\gamma,1)^2\, (c\tau+d)^{1+2\epsilon}\, \theta_\chi\Big(\dfrac{\tau}{N}\Big)\, \theta_{\bar\chi}\Big(\dfrac{\tau}{N}\Big).
$$
Hence
$$
B_{\chi}(\gamma(\tau))\=\dfrac{\upsilon(\gamma,1)^2}{\epsilon_2(\gamma')^{2(1+2\epsilon)}}\, B_{\chi}(\tau).
$$

We use the congruences \eqref{congruencias 1} to calculate $\upsilon(\gamma,1)^2/\epsilon_2(\gamma')^{2(1+2\epsilon)}$.

On the one hand, $cdN(1-a^2)\equiv 0\pmod{3}$ because either $a^2\equiv 1\pmod{3}$, either $a\equiv 0\pmod{3}$, in which case $d\equiv 0\pmod{3}$.

On the other hand, $a\frac{b}{N}-acN\equiv-acN(1+\mathrm{det}(\mu))\equiv 0\pmod{12}$. The first congruence is clear, also is the second modulo 4. For the second congruence modulo 3, either $\mathrm{det}(\mu)\equiv-1\pmod{3}$, either $\mathrm{det}(\mu)\equiv 1\pmod{3}$, in which case $a\equiv d\pmod{3}$, so $ad-bc=1$ implies $bc\equiv 0\pmod{3}$, and thus $c\equiv 0\pmod{3}$ or $N\equiv 0\pmod{3}$.
Hence
\begin{equation}\label{a1}
\epsilon_2(\gamma')^2=\zeta_4^{3cdN(1-a^2)+c_0N(a-1)}.
\end{equation}
We distinguish two cases.

1) If $d$ is odd, then, the congruences \eqref{congruencias 1} and the equation $ad-bc=1$ imply that $ad\equiv 1\pmod{4}$, so the exponent of $\zeta_4$ in \eqref{a1} becomes
$$
c_0N(a-1)\equiv c_0N(d-1)\pmod{4}.
$$
Therefore
$$
\dfrac{\upsilon(\gamma,1)^2}{\epsilon_2(\gamma')^{2(1+2\epsilon)}}=\dfrac{\zeta_4^{d-1}}{\zeta_4^{c_0N(d-1)(1+2\epsilon)}}=\zeta_4^{(d-1)(1-c_0N(1+2\epsilon))}=1.
$$

2) If $d$ is even, then, because of congruences \eqref{congruencias 1}, the exponent of $\zeta_4$ in \eqref{a1} becomes
$$
3cdN+cN(a-1)\equiv cN(3d+a-1)\equiv -cN\pmod{4},
$$
so 
$$
\dfrac{\upsilon(\gamma,1)^2}{\epsilon_2(\gamma')^{2(1+2\epsilon)}}=\dfrac{\zeta_4^{b}}{\zeta_4^{-cN(1+2\epsilon)}}=\zeta_4^{c((1+2\epsilon)N-1)}=(-1)^{\frac{N-v}{2}}.
$$
The second equality can be deduced from the congruences $b\equiv c\equiv 1\pmod{2}$ and $bc\equiv -1\pmod{4}$.

Thus
$$
(B_{\chi}|\mu)(\tau)=(-1)^{\frac{N-v}{2}(t-1)}\, B_{\chi^{\mathrm{det}(\mu)}}(\tau).
$$

We can explicit the action of $W_{M,iN}$ on $A_{\chi}(\tau)^n$ in a similar way. The expression \eqref{expresion} becomes in this case
$$
\upsilon(\gamma',N)^{n}\, (c\tau+d)^{n/2+n\epsilon}\sum_{h_1,\ldots, h_{n}\, \mathrm{mod}\, N}\,\, \prod_{j=1}^n \, \chi(a^{-1}h_j)\, \theta_{N,h_j}^{(\epsilon)}\Big(\dfrac{\tau}{N}\Big),
$$
with $\upsilon(\gamma',N)=\upsilon(\gamma,1).$
Since $\chi(a^{-1})^{n}=1$, following the previous notations \eqref{a}, we have
$$
\theta_\chi\Big(\gamma'\Big(\dfrac{\tau}{N}\Big)\Big)^n=\upsilon(\gamma,1)^{n}\, (c\tau+d)^{n/2+n\epsilon}\, \theta_\chi\Big(\dfrac{\tau}{N}\Big)^n,
$$
so
$$
A_{\chi}(\gamma(\tau))^n=\frac{\upsilon(\gamma,1)^n}{(\epsilon_1(\gamma')\epsilon_2(\gamma'))^{(1+2\epsilon)n}}\,  A_{\chi}(\tau)^n=(-1)^{\frac{(N-v)n}{2}(d-1)}A_{\chi}(\tau)^n
$$
and
$$
(A_{\chi}|\mu)(\tau)^n=(-1)^{\frac{(N-v)n}{2}(t-1)}\, A_{\chi^{\mathrm{det}(\mu)}}(\tau)^n.
$$
\begin{flushright}
$\square$
\end{flushright}

From now on we suppose $N=p>2$ is prime and we denote by $X(p,m)$ the set of characters with conductor $p$ and order $m$ up to complex conjugation.
All characters with fixed prime conductor and fixed order have the same parity; as before $v=1$ if they are even and $v=-1$ if they are odd. 

\begin{teorema}\label{teo final 2} The following sets are orbits for the action of the group $\mathrm{Gal}(H_{M,\O}/H_\O)$ on the field $H_{M,\O}$:
\begin{align*}
(i)\,  &\left\{B_{\chi}(ip)^2\mid\chi\in X(p,m)\right\},\\
(ii)\, &\left\{B_{\chi}(ip)\mid\chi\in X(p,m)\right\} \, \mbox{if $p\equiv v\!\!\pmod 4$},\\
(ii)\, &\left\{A_{\chi}(ip)^{2n}, A_{\bar{\chi}}(ip)^{2n}\mid\chi\in X(p,m)\right\},\\
(iv)\, &\left\{A_{\chi}(ip)^{n}, A_{\bar{\chi}}(ip)^{n}\mid\chi\in X(p,m)\right\} \, \mbox{if $m\equiv 0\!\!\!\pmod 4$ or $p\equiv v\!\!\!\pmod 4$}.
\end{align*}
\end{teorema}

The proof follows from the two lemmas below.

\begin{lema}\label{lema orbita 1} Given $\chi\in X(p,m)$, we have
$$
X(p,m)\=\left\{\chi^\sigma,\, \bar{\chi}^\sigma\, \mid\, \sigma\in(\mathbb{Z}/m\mathbb{Z})^\ast\right\}.
$$ 
\end{lema}

\textbf{Proof.}  The inclusion of the right-hand set into $X(p,m)$ is clear. We should see that given $\chi$ and $\chi'$ in $X(p,m)$, the character $\chi'$ is in the $(\mathbb{Z}/m\mathbb{Z})^\ast$-orbit of $\chi$.

The group $(\mathbb{Z}/p\mathbb{Z})^\ast$ is cyclic; let $h$ be a generator. The groups $\mathrm{Im}(\chi)$ and $\mathrm{Im}(\chi')$ are contained in the group of the m-th roots of unity, which is also cyclic and from which $\chi(h)$ and $\chi'(h)$ are generators. We write $\chi'(h)=\chi(h)^\sigma$ with $\sigma\in(\mathbb{Z}/m\mathbb{Z})^\ast$. 

For $h^{\sigma'}\in(\mathbb{Z}/p\mathbb{Z})^\ast$, we have 
$$
\chi'(h^{\sigma'})=\chi'(h)^{\sigma'}=\chi(h)^{\sigma\sigma'}=\chi(h^{\sigma'})^\sigma,
$$ 
so $\chi'=\chi^\sigma$.
\begin{flushright}
$\square$
\end{flushright}

\begin{lema}\label{lema orbita 2} The following sets equality is satisfied:
\begin{equation}\label{orbita 2}
(\mathbb{Z}/m\mathbb{Z})^\ast=\left\{\pm(t^2+s^2)\pmod{m}:\,  (t^2+p^2s^2,6mp)=1 \right\}.
\end{equation} 
\end{lema}

\textbf{Proof.} Let $u\in\mathbb{Z}$ be coprime with $m$. By Dirichlet's Theorem, there exists a prime number $q\neq 3,p$ such that
$$
q\equiv\left\{\begin{array}{lll}
u\!\!\pmod{4m}  &\quad\mbox{if $u\equiv 1\!\!\pmod{4}$}\\
\!-u\!\!\pmod{4m} &\quad\mbox{if $u\equiv 3\!\!\pmod{4}$}\\
u+m\!\!\pmod{4m} &\quad\mbox{if $u\equiv 0\!\!\!\pmod{2}$ and $u+m\equiv 1\!\!\!\pmod{4}$}\\
\!-(u+m)\!\!\pmod{4m} &\quad\mbox{otherwise}.
\end{array}\right.
$$
In all cases $q\equiv 1\pmod{4}$ and we can write $q=t^2+s^2$ with $t,s\in\mathbb{Z}$. Hence
$$
u\equiv\pm(t^2+s^2)\pmod{m}.
$$  
We want to show $(t^2+p^2s^2,6pm)=1$. Since all the expressions above are symmetric in $t$ and $s$, we can suppose $p\nmid t$. Then $(t^2+p^2s^2,p)=1$. When $p\neq 3$, the integers $s^2$ and $p^2s^2$ are the same modulo 2, and also modulo 3, so $t^2+p^2s^2\equiv q\pmod{6}$. Since $q\neq 2,3$, $(t^2+p^2s^2,6)=1$. (If $p=3$, also $(t^2+p^2s^2,6)=1$ because $p\nmid t$). Since $p\equiv 1\pmod{m}$, $t^2+p^2s^2\equiv\pm u\pmod{m}$. We chose $u$ coprime with $m$, so $(t^2+p^2s^2,6pm)=1$. Therefore $\pm(t^2+s^2)\pmod{m}$ belongs to the set on the right hand side of  \eqref{orbita 2}.
\begin{flushright}
$\square$
\end{flushright}

By Lemmas \ref{lema orbita 1} and \ref{lema orbita 2}, 
$$
X(p,m)\=\left\{\chi^{\det(\mu)}\, \mid\, \mu\in W_{M,ip}\right\}.
$$ 
Then Theorem \ref{teo final 2} follows from Proposition \ref{prop 2}.

\begin{coro} We have
\begin{align*}
&(i)\,\, \mathcal{N}(p,m)^2\in H_\O,\\
&(ii)\,\, \mathcal{N}(p,m)\in H_\O \quad\mbox{if\, $|X(p,m)|\equiv 0\pmod{2}$ \, or\, $p\equiv v\pmod{4}$}.
\end{align*}
\end{coro}

\begin{coro} For all $\chi\in X(p,m)$,
$$
[H_\O(B_\chi(ip)):K]\leq\left\{\begin{array}{ll}
\dfrac{|X(p,m)|(p-v)}{2} &\mbox{if $p\equiv v\pmod{4}$},\\
&\\
|X(p,m)|(p+v) &\mbox{if $p\equiv -v\pmod{4}$}.
\end{array}\right.
$$
\end{coro}
 
\textbf{Proof.} Denoting by $h(\mathcal{O})$ the class number of $\mathcal{O}=\mathbb{Z}[ip]$,
$$
[H_\O:K]\=|\mathrm{Gal}(H_\O/K)|\=|Cl(\mathcal{O})|\=h(\mathcal{O}).
$$
Applying the general formula for the class number of an imaginary quadratic order (see \cite{C}), we have
$$
h(\mathcal{O})=\left\{\begin{array}{ll}\dfrac{p-1}{2}&\mbox{if $p\equiv 1\pmod{4}$}\\
&\\
\dfrac{p+1}{2}&\mbox{if $p\equiv 3\pmod{4}$}.
\end{array}\right.
$$   
We deduce from the statements (i) and (ii) of Theorem \ref{teo final 2}
$$
[H_\O(B_\chi(ip)):H_\O]\leq\left\{\begin{array}{ll}
2|X(p,m)| &\mbox{if $p\equiv -v\pmod{4}$}\\
&\\
|X(p,m)| &\mbox{if $p\equiv v\pmod{4}$}.
\end{array}\right.
$$
\begin{flushright}
$\square$
\end{flushright}



\begin{teorema}
There is a constant $c>0$ such that for all non-quadratic $\chi$ with prime conductor $p=2l+1$, where $l$ is prime, satisfying $p>c$, we have $\theta_\chi(i)\neq 0$.
\end{teorema}

\textbf{Proof.} Louboutin proved in \cite{Lo} that there is a constant $c>0$ such that $\theta_\chi(i)\neq 0$ for at least $cp/\log(p)$ characters of the $(p-1)/2$ odd characters with conductor $p$ and of the $(p-1)/2$ even ones. When $p=2l+1$, there is one odd character having order $2$, $(p-3)/2$ odd characters having order $2l$, $(p-1)/2$ even characters having order $l$ and the trivial (even) character. By Theorem \ref{teo final 2}, if $B_\chi(ip)\neq 0$ for some $\chi\in X(p,m)$, then $B_\chi(ip)\neq 0$ for all $\chi\in X(p,m)$. Thus $\theta_\chi(i)\neq 0$ for all non-quadratic characters with conductor $p$ satisfying $\log(p)/p<c$.      
\begin{flushright}
$\square$
\end{flushright}

\begin{rk} For odd but maybe not prime $N$, Theorem \ref{teo final 2} does not apply, but we have
\begin{equation}\label{(i)}
\prod_{\chi\in X(N,m)} (X-B_\chi(iN)^2)\, \in H_\O[X].
\end{equation}
If $N|X(N,m)|\equiv\sum_{\chi\in X(N,m)}\chi(-1)\pmod{4}$, then the square in \eqref{(i)} is not necessary.
\end{rk}

\textbf{Acknowledgements.} This work is part of my PhD thesis. I wish to express my gratitude to Don Zagier for his precious advice in discussing mathematics and to Pilar Bayer for her careful reading of this paper.


\begin{thebibliography}{9}

\bibitem[CZ13]{C-Z}
Cohen, H; Zagier, D.: \emph{Vanishing and non-vanishing theta values}. Annales math�matiques du Qu�bec 37 (2013), 45-61 (special issue dedicated to Professor Paulo Ribenboim).

\bibitem[Cox89]{C}
Cox, D.: \emph{Primes of the form $x^2+ny^2$}. John Wiley \& Sons, 1989. 

\bibitem[Fra35]{Fra}
Franz, W.: \emph{Die Teilwert der Weberschen Tau-Funktion}. J. reine angew. Math. 173 (1935) 60-64.

\bibitem[Gee00]{G1}
Gee, A.: \emph{Class fields by Shimura reciprocity}. Thesis, University of Amsterdam (2000).

\bibitem[GS98]{GS}
Gee, A.; Stevenhagen, P.: \emph{Generating class fields using Shimura reciprocity}. Algorithmic Number Theory (J. P. Buhler, ed.). Springer LNCS 1423 (1998) 441-453.

\bibitem[Iwa97]{I}
Iwaniec, H.: \emph{Topics in Classical Automorphic Forms}, American Mathematical Soc, 1997.

\bibitem[Lan87]{La2}
Lang, S.: \emph{Elliptic functions}, 2nd edition, Springer-Verlag, Berlin-Heidelberg-New York, 1987.

\bibitem[Lou99]{Lo}
Louboutin, S.: \emph{Sur le calcul num�rique des constantes des �quations fonctionnelles des fonctions L associ�es aux caract�res impairs}. C.R. Acad. Sci. Paris 329 (1999) 347-350. 


\bibitem[Mey57]{Me}
Meyer, C.: \emph{�ber einige Anwendungen Dedekindscher Summen}. J.reine angew.Math. 198 (1957) 143-203.

\bibitem[Mum83]{Mu}
Mumford, D.: \emph{Tata Lectures on Theta I}. Birkh�user, 1983.


\bibitem[Shi71]{Sh}
Shimura, G.: \emph{Introduction to the Arithmetic Theory of Automorphic Functions}. Iwanami Shoten and Princeton University Press, 1971.

\bibitem[Ste00]{St}
Stevenhagen, P.: \emph{Hilbert's 12th problem, complex multiplication and Shimura reciprocity}. Adv. studies in pure mathematics, Math. Soc. Japan (2000) 1-16.

\end{thebibliography}
\end{document}